\documentclass[letterpaper, 10 pt, conference]{ieeeconf}
\usepackage{graphics}
\usepackage{epsfig}
\usepackage{cite}
\usepackage{breakurl}
\usepackage{amsmath}
\usepackage{times}

\usepackage{amsthm}
\usepackage{bm}
\usepackage{amssymb}
\usepackage{amsfonts}
\usepackage{bm}
\usepackage{amssymb}
\usepackage{amsfonts}
\usepackage{graphicx}
\usepackage{makecell,multirow,diagbox}
\usepackage{footnote}
\usepackage[normalem]{ulem}
\usepackage{color}
\usepackage{dsfont}
\usepackage[ruled,linesnumbered]{algorithm2e}
\usepackage[dvipsnames]{xcolor}
\usepackage{booktabs, tabularx}
\IEEEoverridecommandlockouts
\newtheorem{customthm}{Theorem}
\newtheorem{customlem}{Lemma}
\newtheorem{customass}{Assumption}

\newtheorem{customdef}{Definition}
\newtheorem{customrem}{Remark}
\allowdisplaybreaks[0]

\begin{document}

\title{{\LARGE \textbf{Joint Time and Energy-Optimal Control of Connected Automated
Vehicles at Signal-Free Intersections with Speed-Dependent Safety Guarantees}}}
\author{Yue Zhang, and Christos G. Cassandras \thanks{Supported in part by NSF under
grants ECCS-1509084, DMS-1664644, CNS-1645681, by AFOSR under grant
FA9550-15-1-0471, by ARPA-E's NEXTCAR program under grant {DE-AR}0000796 and by
the MathWorks.} \thanks{Y. Zhang and C.G. Cassandras are with the Division of
Systems Engineering and Center for Information and Systems Engineering, Boston
University, Boston, MA 02215 USA (e-mail: joycez@bu.edu; cgc@bu.edu).} }
\maketitle


\begin{abstract}
We extend earlier work establishing a framework for optimally
controlling Connected Automated Vehicles (CAVs) crossing a signal free
intersection by \emph{jointly} optimizing energy and travel time. We derive
explicit optimal control solutions in a decentralized manner that guarantee both a speed-dependent
rear-end safety constraint and a time-dependent lateral collision constraint,
in addition to lower/upper bounds on speed and acceleration. Extensive
simulation examples are included to illustrate this framework.
\end{abstract}


\section{Introduction}

\label{intro}

Traffic control at intersections is one of the major challenges in
transportation systems as intersections account for a large fraction of
accidents and of the overall system congestion. To date, traffic light
control is the prevailing method for coordinating conflicting traffic flows
through an intersection. Recent technological developments include designing
online adaptive traffic light control as in \cite{fleck2016adaptive}.
However, aside from the obvious infrastructure cost of traffic lights, the
efficiency and safety offered by such signaling methods can be significantly
improved through new approaches capable of enabling smoother traffic flow
while ensuring safety.

Connected Automated Vehicles (CAVs) provide the most intriguing opportunity
for improving traffic conditions in a transportation network. One of the very early efforts
was proposed in \cite{Athans1969} and \cite{Levine1966} where the a linear
optimal regular is introduced to control a single string of vehicles for the
merging problem. More recently, 
Dresner and Stone \cite{Dresner2004} proposed a reservation-based scheme for
automated intersection management. Since then, numerous research efforts
have explored efficient and safe control strategies, e.g., \cite%
{Dresner2008, DeLaFortelle2010, Huang2012}. Some of the efforts focused on
minimizing travel delays with safety guarantees \cite%
{Li2006,Yan2009,Zhu2015,Zohdy2012, Lee2012, Miculescu2014}. 
Lee and Park \cite{Lee2012} aimed at minimizing the overlap between vehicle
positions. Miculescu and Karaman \cite{Miculescu2014} have studied
intersections as polling systems and determined a sequence of times assigned
to vehicles on each road. 
Reducing energy consumption is another desired objective which has been
considered in recent literature \cite{gilbert1976vehicle, hooker1988optimal,
hellstrom2010design, li2012minimum}. Hellstrom \cite{hellstrom2010design}
proposed an energy-optimal control algorithm for heavy diesel trucks by
utilizing road topography information. 
A detailed discussion of recent advances in this area can be found in \cite%
{rios2017survey}.

The contribution of this paper consists of extending the optimal control
framework in \cite{malikopoulos2018decentralized}. First, instead of solving
a throughput maximization problem followed by an energy minimization problem
for each CAV, here we formulate a problem in which each CAV seeks to \emph{%
jointly minimize} both its travel time through a specified \textit{Control
Zone} (CZ) and \textit{Merging Zone} (MZ) and its energy consumption. This
allows us to readily quantify the tradeoff between these two criteria (see
also \cite{arxiv1} where left/right turns are included along with a
passenger comfort metric). Second, unlike \cite{malikopoulos2018decentralized}, \cite{arxiv1} where we first resolve
possible collisions in the MZ and then apply optimal control over the CZ,
here we relax the constant speed assumption inside the MZ and handle lateral
collision avoidance as additional state constraints; this provides
flexibility in controlling CAVs within the MZ. Third, unlike \cite%
{malikopoulos2018decentralized}, \cite{arxiv1} where we limit ourselves to
a distance-dependent rear-end safety constraint, here we include a
speed-dependent rear-end safety constraint, which better captures the
relationship between two consecutive vehicles traveling on the same road.
Our analysis includes the derivation of several structural properties of an
optimal control solution and it allows us to determine whether an optimal
control solution for each CAV is feasible at the time it enters the CZ.

The paper is structured as follows. In Section II, we review the model in 
\cite{malikopoulos2018decentralized} and derive the conditions that
guarantee safety constraints for each CAV. In Section III, we formulate a
decentralized optimal control problem for each CAV that jointly minimizes
its travel time and energy consumption throughout the CZ and the MZ, prove
structural properties of optimal trajectories, and derive an explicit
solution for it. Simulation results are given in Section IV showing
constrained optimal trajectories with different safety constraints becoming
active. Concluding remarks are given in Section V.


\section{The Intersection Model}

\label{model}

We begin with a brief review of the model introduced in \cite%
{ZhangMalikopoulosCassandras2016} and fully developed in \cite%
{malikopoulos2018decentralized}. We consider an intersection (Fig. \ref%
{fig:intersection}) where the region at its center, assumed to be a square
of side $S$, is called \emph{Merging Zone} (MZ) and defines the area of
potential lateral CAV collisions. The intersection has a \emph{Control Zone}
(CZ) and the road segment from the CZ entry to the CZ exit (i.e., the MZ
entry) is referred to as a CZ segment whose length $L>S$ is assumed to be
the same for all entry points to a given CZ. Extensions to asymmetric CZ
segments are possible and considered in \cite{Zhang2018sequence}.

\begin{figure}[ptb]
\centering
\includegraphics[width = 0.7 \columnwidth]{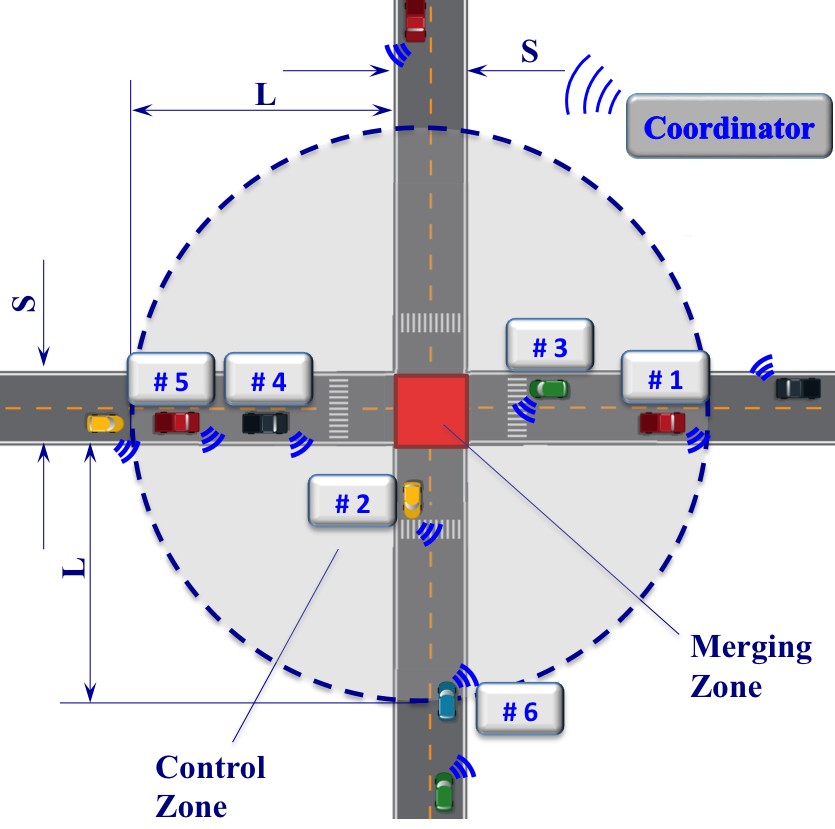}
\caption{Connected Automated Vehicles crossing an urban intersection.}
\label{fig:intersection}
\end{figure}

We assume the existence of a \textquotedblleft coordinator" whose task is to
handle the information exchanges between CAVs, while each CAV maintains its
own control autonomy. Let $N(t)\in \mathbb{N}$ be the cumulative number of
CAVs which have entered the CZ by time $t$ and formed a queue that
designates the crossing sequence in which these CAVs will enter the MZ.
There is a number of ways to manage such a queue. In \cite%
{malikopoulos2018decentralized} a strict First-In-First-Out (FIFO) crossing
sequence is assumed, that is, when a CAV reaches the CZ, the coordinator
assigns it an integer value $i=N(t)+1$. This is relaxed in \cite%
{Zhang2018sequence} to allow for dynamically resequencing CAVs as each new
one arrives, hence maximizing throughput. If two or more CAVs enter a CZ at
the same time, then the corresponding coordinator selects randomly the first
one to be assigned the value $N(t)+1$.

For simplicity, we assume that each CAV is governed by second order dynamics:%
\begin{equation}
\dot{p}_{i}=v_{i}(t)\text{, }~p_{i}(t_{i}^{0})=0\text{; }~\dot{v}%
_{i}=u_{i}(t)\text{, }v_{i}(t_{i}^{0})\text{ given}   \label{eq:model2}
\end{equation}
where $p_{i}(t)\in\mathcal{P}_{i}$, $v_{i}(t)\in\mathcal{V}_{i}$, and $%
u_{i}(t)\in\mathcal{U}_{i}$ denote the position, i.e., travel distance since
the entry of the CZ, speed and acceleration/deceleration (control input) of
each CAV $i$. The sets $\mathcal{P}_{i}$, $\mathcal{V}_{i}$ and $\mathcal{U}%
_{i}$ are complete and totally bounded subsets of $\mathbb{R}$. These
dynamics are in force over an interval $[t_{i}^{0},t_{i}^{f}]$, where $%
t_{i}^{0}$ and $t_{i}^{f}$ are the times that the vehicle $i$ enters the CZ
and exits the MZ respectively. To ensure that the control input and vehicle
speed are within a given admissible range, the following constraints are
imposed: 
\begin{equation}
\begin{split}
u_{i,min} & \leq u_{i}(t)\leq u_{i,max},\quad\text{and} \\
0 & \leq v_{min}\leq v_{i}(t)\leq v_{max},\quad\forall t\in\lbrack
t_{i}^{0},t_{i}^{f}].
\end{split}
\label{speed_accel constraints}
\end{equation}


As part of safety considerations, we impose the following assumptions:

\begin{customass}
CAVs follow the crossing sequence established by the coordinator and no
overtaking, reversing directions, lane-changing, or turns are allowed. 
\label{ass:infra}
\end{customass}

\begin{customass}
Each vehicle has proximity sensors and can observe and/or estimate local
information that can be shared with other vehicles. 
\label{ass:sensor}
\end{customass}

\begin{customass}
For each CAV, the speed constraints in \eqref{speed_accel constraints} and
the rear-end safety constraint in \eqref{rearend} are not active at $%
t_{i}^{0}$. \label{ass:feas}
\end{customass}

If the last assumption is violated, any optimal control solution is
obviously infeasible and we must resort to control actions that simply
attempt to satisfy these constraints as promptly as possible; alternatively,
we may impose a Feasibility Enforcement Zone (FEZ) that precedes the CZ as
described in \cite{Zhang2017feasibility}.

\begin{customdef}
Depending on its physical location inside the CZ, CAV $i-1\in\mathcal{N}(t)$
belongs to only one of the following four subsets of $\mathcal{N}(t)$ with
respect to CAV $i$: 1) $\mathcal{R}_{i}(t)$ contains all CAVs traveling on
the same road as $i$ and towards the same direction but on different lanes,
2) $\mathcal{L}_{i}(t)$ contains all CAVs traveling on the same road and
lane as vehicle $i$ (e.g., $\mathcal{L}_{5}(t)$ contains CAV \#4 in Fig. \ref%
{fig:intersection}), 3) $\mathcal{C}_{i}(t)$ contains all CAVs traveling on
different roads from $i$ and having destinations that can cause collision at
the MZ (e.g., $\mathcal{C}_{6}(t)$ contains CAV \#5 in Fig. \ref%
{fig:intersection}), and 4) $\mathcal{O}_{i}(t)$ contains all CAVs traveling
on the same road as $i$ and opposite destinations that cannot, however,
cause collision at the MZ (e.g., $\mathcal{O}_{4}(t)$ contains CAV \#3 in
Fig. \ref{fig:intersection}). \label{def:2}
\end{customdef}

A rear-end collision may occur only if some CAV $z\neq i$ belongs to $%
\mathcal{L}_{i}(t)$. To ensure the absence of any rear-end collision
throughout the CZ and MZ, instead of using the distance-dependent rear-end
safety constraint as in \cite{malikopoulos2018decentralized},\cite%
{Zhang2017turns}, we impose a speed-dependent rear-end safety constraint 
\begin{equation}
\begin{aligned} s_{i}(t)& =p_{k}(t)-p_{i}(t)\geq \varphi v_i(t) +
\delta_0,\\ \forall t & \in\lbrack t_{i}^{0},t_{i}^{f}],
~k=\max_{z}\{z\in\mathcal{L}_{i}(t)\}\label{rearend}\end{aligned}
\end{equation}%
that specifies a minimum safe headway, i.e., a gap that is a function of $%
v_{i}(t)$, where $k$ is the CAV physically ahead of $i$, $\varphi $ is the
reaction time and $\delta _{0}$ is the minimal standstill inter-vehicle distance. Note that in \cite{malikopoulos2018decentralized} where we use a distance-dependent safety
constraint $p_{k}(t)-p_{i}(t)\geq \delta $, $\delta $ denotes an
inter-vehicle distance while vehicles are moving.

As in \cite{malikopoulos2018decentralized}, we consider a First-In-First-Out
(FIFO) ordering structure by imposing the following condition: 
\begin{equation}
t_{i}^{f}\geq t_{i-1}^{f}  \label{fifo}
\end{equation}

A lateral collision involving CAV $i$ may occur only if some CAV $z\neq i$
belongs to $\mathcal{C}_{i}(t)$. Letting $t_{i}^{m}$ denote the time when a
CAV enters the MZ, this leads to the following definition:

\begin{customdef}
For each CAV $i\in \mathcal{N}(t)$, the set $\Gamma_{i}$ includes all time
instants when a lateral collision involving CAV $i$ is possible: 
\begin{equation}
\Gamma _{i}\triangleq \Big\{t~|~t\in \lbrack t_{i}^{m},t_{i}^{f}]\Big\}. 
\notag
\end{equation}
\end{customdef}

Consequently, to avoid a lateral collision for any two vehicles $i,z\in 
\mathcal{N}(t)$ on different roads, the following constraint should hold 
\begin{equation}
\Gamma _{i}\cap \Gamma _{z}=\varnothing ,\text{ \ }\forall t\in \lbrack
t_{i}^{m},t_{i}^{f}]\text{, \ }z\in \mathcal{C}_{i}(t).  \label{eq:lateral}
\end{equation}%

Combining \eqref{eq:lateral} with the FIFO constraint \eqref{fifo}, we can
easily derive the lateral collision constraint 
\begin{equation}
t_{i}^{m}\geq t_{c}^{f},~c=\max_{z}\{z\in \mathcal{C}_{i}(t)\}
\label{lateral}
\end{equation}
Since $t_{i}^{m}$ is the time that CAV $i$ reaches the end of the CZ, the
constraint \eqref{lateral} is equivalent to the following position-dependent interior-point constraint 
\begin{equation}
p_{i}(t_{c}^{f})\leq L, ~c=\max_{z}\{z\in \mathcal{C}_{i}(t)\}.   \label{lateral2}
\end{equation}%
For CAVs that belong to $\mathcal{O}_{i}(t)$ and $\mathcal{R}_{i}(t)$, no
collision would occur throughout the CZ and the MZ. 

To ensure that CAV $k$ will not collide with CAV $i$ after $k$ exits the MZ while $i$ is still traveling inside the MZ, we apply the following assumption.
\begin{customass}
The speed of CAV $k$ will remain constant after exiting the MZ for $t \in [t_{k}^f, t_i^f]$, $k=\max_{z}\{z\in\mathcal{L}_{i}(t)$. 
\label{ass:constant} 
\end{customass} 

Combining \eqref{speed_accel constraints}, \eqref{rearend}, \eqref{fifo},
the lower bound on the time when CAV $i$ can leave the MZ is given by 
\begin{equation}
t_{i}^{f}\geq \max \{t_{i}^{L},t_{k}^{f}+\frac{\varphi v_{i}^{f}+\delta _{0}%
}{v_{k}^{f}},t_{o}^{f}\},  \label{tif_lw}
\end{equation}%
where $k=\max_{z}\{z\in \mathcal{L}_{i}(t)\}$, $o=\max_{z}\{z\in \mathcal{O}%
_{i}(t)\cup \mathcal{R}_{i}(t)\}$, and $t_{i}^{L}=t_{i}^{1}\mathds{1}%
_{v_{i}^{f}=v_{max}}+t_{i}^{2}(1-\mathds{1}_{v_{i}^{f}=v_{max}})$ is the
lower bound constrained by the speed and control constraints and $t_{i}^{1}$%
, $t_{i}^{2}$ were derived in \cite{ZhangMalikopoulosCassandras2016}:  
\begin{align}
t_{i}^{1}& =t_{i}^{0}+\frac{L+S}{v_{max}}+\frac{(v_{max}-v_{i}^{0})^{2}}{%
2u_{i,max}v_{max}},  \notag \\
t_{i}^{2}& =t_{i}^{0}+\frac{[2(L+S)u_{i,max}+(v_{i}^{0})^{2}]^{1/2}-v_{i}^{0}%
}{u_{i,max}}.  \label{t_lw}
\end{align}

\begin{customthm}
If CAV $z$, $z\in\{2,\ldots,i\}$ satisfies \eqref{rearend}, \eqref{lateral},
and \eqref{tif_lw}, then, with respect to any CAV $j$, $j<i$, CAV $i$ satisfies

\begin{itemize}
\item $p_{j}(t)-p_{i}(t)\geq\varphi v_{i}(t)+\delta_{0}$, $t\in\lbrack
t_{i}^{0},t_{i}^{f}]$ if $j\in\mathcal{L}_{i}(t)$ [no rear-end collision],
\item $t_{i}^{m}\geq t_{j}^{f}$ if $j\in\mathcal{C}_{i}(t)$ [no lateral
collision inside the MZ],
\item $t_{i}^{f}\geq t_{j}^{f}$ [crossing order preservation].

\end{itemize}

\label{Theorem1}
\end{customthm}

\emph{Proof. } If $\max \{t_{i}^{L},t_{k}^{f}+\frac{\varphi v_{i}^{f}+\delta
_{0}}{v_{k}^{f}},t_{o}^{f}\}=t_{i}^{L}$, then $t_{i}^{f}\geq t_{i}^{L}$
ensures that $t_{i}^{f}$ is feasible since it depends only on the control
and state constraints \eqref{speed_accel constraints}. Based on the definition of CAV $k$, $c$ and $o$, CAV $i$ will not generate any collision with $j$ if $j \in \{k, c, o\}$. When $j\neq
k,c,o$, there are three cases to consider as follows. 

(1) When $j\in \mathcal{L}_{i}(t)$. In this case, a rear-end collision is
possible. Since $k,j\in \mathcal{L}_{i}(t)$ and $j<k$, CAVs $k$ and $j$ are
traveling on the same lane towards the same direction as $i$, we have $j<k<i$%
, and $p_{i}(t)+\varphi v_{i}(t)+\delta _{0}\leq p_{k}(t)<p_{k}(t)+\varphi
v_{k}(t)+\delta _{0}\leq p_{j}(t)$. The ordering is therefore implicitly guaranteed. 

(2) When $j\in \mathcal{C}_{i}(t)$. In this case, only a lateral collision inside
the MZ is possible. Since $c,j\in \mathcal{C}_{i}(t)$ and $j<c$,
there are two subcases to consider: (i) $j\in \mathcal{L}_{c}(t)$, (ii) $%
j\in \mathcal{O}_{c}(t)\cup \mathcal{R}_{c}(t)$. When $j\in \mathcal{L}%
_{c}(t)$, the rear-end safety constraint \eqref{rearend} leads to $t_{c}^{f}>t_{j}^{f}$.
Due to the lateral collision constraint \eqref{lateral}, we have $t_{i}^{m}\geq t_{c}^{f}$. Combining with $t_{c}^{f}>t_{j}^{f}$, we have $t_{i}^{m}>t_{j}^{f}$. When $j\in \mathcal{O}_{c}(t)\cup \mathcal{R}_{c}(t)$, denote CAV $j_o = \max_{z}\{z\in \mathcal{O}_{c}(t)\cup \mathcal{R}_{c}(t)\}$. If $j = j_o$, then we have $t_{c}^{f}\geq t_{j}^{f}$ following the third term in \eqref{tif_lw}; if $j < j_o$, CAV $j$ must exit the MZ earlier than CAV $j_o$, and hence, $t_{c}^{f}\geq t_{j_o}^{f} \geq t_j^f$. Combining with $t_{i}^{m}\geq t_{c}^{f}$, we have $t_{i}^{m}\geq t_{j}^{f}$. The ordering is therefore implicitly guaranteed.

(3) When $j\in \mathcal{O}_{i}(t)\cup \mathcal{R}_{i}(t)$. In this case, no
collision between $i$ and $j$ could
occur according to the definition. We only need to ensure the CAV ordering.
Since $o,j\in \mathcal{O}_{i}(t)\cup \mathcal{R}_{i}(t)$ and $j<o$, there
are two subcases to consider: (i) $j\in \mathcal{L}_{o}(t)$, (ii) $j\in 
\mathcal{O}_{o}(t)\cup \mathcal{R}_{o}(t)$. When $j\in \mathcal{L}_{o}(t)$,
the rear-end safety constraint \eqref{rearend} ensures the CAV ordering, i.e., $t_{o}^{f}>t_{j}^{f}$. Following the third term in \eqref{tif_lw}, we have $t_{i}^{f}\geq t_{o}^{f}$. Combining with $t_{o}^{f}>t_{j}^{f}$, we can obtain $t_{i}^{f}>t_{j}^{f}$, which satisfies the ordering constraint. When $j\in \mathcal{O}_{o}(t)\cup \mathcal{R}_{o}(t)$, we
have $t_{i}^{f}\geq t_{o}^{f}\geq t_{j}^{f}$ following the third term in \eqref{tif_lw}, which completes the proof.

\hfill$\blacksquare$

Corresponding to the lower bound of terminal time $t_{i}^{L}$, there also
exists the upper bound $t_{i}^{U}$: 
\begin{equation}
t_{i}^{U}=t_{i}^{3}\mathds{1}_{v_{i}^{f}=v_{min}}+t_{i}^{4}(1-\mathds{1}%
_{v_{i}^{f}=v_{min}})  \label{t_up}
\end{equation}%
where $v_{i}(t_{i}^{f})=\sqrt{2(L+S)u_{min}+(v_{i}^{0})^{2}},$ and $%
t_{i}^{3}=t_{i}^{0}+\frac{L+S}{v_{min}}+\frac{(v_{min}-v_{i}^{0})^{2}}{%
2u_{min}v_{min}}$ and $t_{i}^{4}=t_{i}^{0}+\frac{v_{i}(t_{i}^{f})-v_{i}^{0}}{%
u_{min}}$ are derived in a similar way as $t_{i}^{1}$ and $t_{i}^{2}$ in %
\eqref{t_lw} respectively (see \cite{ZhangMalikopoulosCassandras2016}).
Based on \eqref{t_up}, the following upper bound constraint applies: 
\begin{equation}
t_{i}^{f}\leq t_{i}^{U}  \label{tm_up_constraint}
\end{equation}


\section{Optimal Control of CAVs in the CZ \label{CZ}}

The objective of each CAV inside the CZ and MZ, i.e., over $%
[t_{i}^{0},t_{i}^{f}]$, is to derive an optimal acceleration/deceleration
which minimizes a convex combination of its travel time and energy
consumption. Since the coordinator is not involved in any decision making
process regarding vehicle control, we can formulate a tractable
decentralized problem, that can be solved on line by each CAV, as follows: 
\begin{equation}
\begin{aligned} &\min_{u_{i}\in U_{i}}\int_{t_{i}^{0}}^{t_{i}^{f}}[\gamma+
\frac{1}{2}u_{i}^{2}(t)]~dt \\ \text{s.t.}:&~\eqref{eq:model2},
(\ref{speed_accel constraints}),\eqref{rearend}, \eqref{lateral2},
\eqref{tif_lw}, (\ref{tm_up_constraint}), \\ & ~p_i(t_i^0)=0,
p_{i}(t_{i}^{f})=L + S \\ & \text{ and given }t_i^0, v_i(t_i^0),
\label{eq:decentral} \end{aligned}
\end{equation}%
where $\gamma $ is a normalized weight associated with the importance of
travel time relative to energy. The constraints consist of the vehicle
dynamics \eqref{eq:model2}, state and control constraints \eqref{speed_accel
constraints}, the speed-dependent rear-end safety constraint \eqref{rearend}%
, the time-dependent lateral constraint \eqref{lateral}, and the lower and
upper bounds of the terminal time $t_{i}^{f}$ \eqref{tif_lw} and %
\eqref{tm_up_constraint}. Unlike the problem considered in \cite%
{malikopoulos2018decentralized} where the terminal time was obtained a
priori to optimize travel times, here the optimal travel time is part of the
problem solution. An additional difference is that the optimization horizon
here covers both CZ and MZ, instead of the CZ only.

\subsection{Problem Decomposition \label{CZ_problem_decomposition}}

In order to efficiently obtain an analytical solution on line, we proceed
with the following step-wise approach (Algorithm \ref{algo:optimization}).
We start with the unconstrained problem $P_{0}$ by relaxing all constraints
in \eqref{eq:decentral} except the dynamics \eqref{eq:model2}. After solving 
$P_{0}$, we obtain the terminal time $t_{i}^{f}$. The second step is to
check whether $t_{i}^{f}$ satisfies both \eqref{tif_lw} and %
\eqref{tm_up_constraint}; if not, we formulate problem $P_1$ by constraining $t_{i}^{f}$ to either the lower
bound \eqref{tif_lw} or the upper bound \eqref{tm_up_constraint} and
re-solve the problem. Then, we proceed with checking the speed, control and
safety constraints \eqref{speed_accel constraints}, \eqref{rearend} and %
\eqref{lateral2} and deal with any violated constraints one by one until they
are all satisfied. Note that if $t_{i}^{L}>t_{i}^{U}$, i.e., the lower bound
on $t_{i}^{f}$ is higher than its upper bound, the problem is obviously
infeasible.

\begin{algorithm}
  formulate an unconstrained problem $P_0$ by relaxing all the constraints in \eqref{eq:decentral} except the dynamics \eqref{eq:model2}\;
  solve $P_0$ and obtain the optimal solution $u_i(t)$\;
  \lIf{$t_i^{f}$ violates \eqref{tif_lw} or \eqref{tm_up_constraint}}{
  formulate $P_1$ by setting $t_i^f$ to either the lower bound \eqref{tif_lw} or the upper bound \eqref{tm_up_constraint};
  solve $P_1$ and obtain a new optimal solution $u_i(t)$}
  \Else{go to step \ref{stepA}}
  set the index of iteration $r := 2$\;
  \Repeat {\eqref{speed_accel constraints}, \eqref{rearend}, and \eqref{lateral2} are all satisfied}{
      check if $u_i(t)$ satisfies \eqref{speed_accel constraints}, \eqref{rearend}, and \eqref{lateral2}\;\label{stepA} 
      \lIf{any of the constraints in \eqref{speed_accel constraints}, \eqref{rearend}, and \eqref{lateral2} is violated}	{
      formulate $P_{r}$ by adding the violated constraint to $P_{r-1}$($P_{1} = P_0$ if $P_1$ does not exist); solve $P_r$ and obtain a new optimal solution $u_i(t)$; 
      $r := r + 1$}
    }
    obtain $u_i(t) = u_i^*(t)$ as the optimal solution\;
  \caption{A step-wise constrained optimization approach}
  \label{algo:optimization}
\end{algorithm}

\subsection{Analytical Solution \label{CZ_solution}}

Given the objective function of the unconstrained problem $P_{0}$, the
Hamiltonian is 
\begin{equation}
\begin{aligned} H_i(p_i, v_i, u_i, \lambda_i, t) = \gamma + \frac{1}{2}
u_i^2(t) + \lambda_i^p v_i(t) + \lambda_i^vu_i(t) \end{aligned}
\label{hamil}
\end{equation}%
and the Lagrangian with constraints directly adjoined is 
\begin{equation}
\begin{aligned} L_i(p_i, v_i, u_i, \lambda_i, \mu_i, \nu_i, t) &= H_i(p_i,
v_i, u_i, \lambda_i, t) + \mu_i g_i(u_i, t) \\ & + \nu_i h_i(p_i, v_i, t) + \zeta_i q_i(p_i(t_1), t_1)
\end{aligned}  \label{lag}
\end{equation}%
where (omitting time arguments for simplicity) $\lambda _{i}=[\lambda
_{i}^{p},\lambda _{i}^{v}]^{T}\in \mathbb{R}^{2}$ is the costate vector, $%
g_{i}(u_{i},t)\leq 0$ and $h_{i}(p_{i},v_{i},t)\leq 0$ represent the control
and state constraints respectively, $q_i(p_i(t_1), t_1) \leq 0$ represents the position-dependent interior-point constraint at $t_1$, and $\mu_{i}=[\mu _{i}^{a},\mu_{i}^{b}]^{T}\in \mathbb{R}^{2}$, $\nu_{i}=[\nu_{i}^{c},\nu_{i}^{d},\nu_{i}^{s}]^{T}\in \mathbb{R}^{3}$, $\zeta_i$ are Lagrange multipliers with 
\begin{equation}
\mu_{i}^{a}=\left\{ 
\begin{array}{ll}
>0, & \mbox{$u_{i}(t) - u_{max} =0$}, \\ 
=0, & \mbox{$u_{i}(t) - u_{max} <0$},%
\end{array}%
\right.   \label{eq:17a}
\end{equation}%
\begin{equation}
\mu_{i}^{b}=\left\{ 
\begin{array}{ll}
>0, & \mbox{$u_{min} - u_{i}(t) =0$}, \\ 
=0, & \mbox{$u_{min} - u_{i}(t)<0$},%
\end{array}%
\right.   \label{eq:17b}
\end{equation}%
\begin{equation}
\nu_{i}^{c}=\left\{ 
\begin{array}{ll}
>0, & \mbox{$v_{i}(t) - v_{max} =0$}, \\ 
=0, & \mbox{$v_{i}(t) - v_{max}<0$},%
\end{array}%
\right.   \label{eq:17c}
\end{equation}%
\begin{equation}
\nu_{i}^{d}=\left\{ 
\begin{array}{ll}
>0, & \mbox{$v_{min} - v_{i}(t)=0$}, \\ 
=0, & \mbox{$v_{min} - v_{i}(t)<0$}.%
\end{array}%
\right.   \label{eq:17d}
\end{equation}%
\begin{equation}
\nu_{i}^{s}=\left\{ 
\begin{array}{ll}
>0, & \mbox{$p_i(t) + \varphi v_i(t) + \delta_0 - p_k(t) = 0$}, \\ 
=0, & \mbox{$p_i(t) + \varphi v_i(t) +  \delta_0 - p_k(t) <0$},%
\end{array}%
\right.   \label{eq:17e}
\end{equation}
\begin{equation}
\zeta_{i}=\left\{ 
\begin{array}{ll}
>0, & \mbox{$p_i(t_c^f) - L = 0$}, \\ 
=0, & \mbox{$p_i(t_c^f) - L <0$},%
\end{array}%
\right.   \label{eq:17f}
\end{equation} 
where $k=\max_{z}\{z\in \mathcal{L}_{i}(t)\}$, $c=\max_{z}\{z\in \mathcal{C}_{i}(t)\}$ and their trajectories including the terminal time
are known to $i$ through the coordinator (or through on-board sensors).

The Euler-Lagrange equations become 
\begin{equation}
\dot{\lambda}_{i}^{p}(t)=-\frac{\partial L_{i}}{\partial p_{i}}=\left\{ 
\begin{array}{ll}
0, & \mbox{$p_i(t) + \varphi v_i(t) + \delta_0 - p_k(t) < 0$}, \\ 
-\nu _{i}^{s}, & \mbox{$p_i(t) + \varphi v_i(t) + \delta_0 - p_k(t) = 0$},%
\end{array}%
\right.   \label{eq:lambda_p}
\end{equation}%
and 
\begin{equation}
\dot{\lambda}_{i}^{v}(t)=-\frac{\partial L_{i}}{\partial v_{i}}=\left\{ 
\begin{array}{ll}
-\lambda _{i}^{p}(t), & \mbox{$v_{i}(t) - v_{max} <0$}~\text{and} \\ 
& \mbox{$v_{min} - v_{i}(t)<0$}, \\ 
-\lambda _{i}^{p}(t)-\nu _{i}^{c}, & \mbox{$v_{i}(t) - v_{max} =0$}, \\ 
-\lambda _{i}^{p}(t)+\nu _{i}^{d}, & \mbox{$v_{min} - v_{i}(t)=0$}.%
\end{array}%
\right.   \label{eq:lambda_v}
\end{equation}%
\textbf{Terminal conditions}. (i) When $t_{i}^{f}$ is free, we have the
following transversality conditions 
\begin{equation}
\lambda _{i}^{v}(t_{i}^{f})=0,\text{ \ \ }H_{i}(t_{i}^{f})=0
\label{eq:costate_trans}
\end{equation}%
(ii) When $t_{i}^{f}$ is constrained by CAV $k$, we set $\psi
_{i}(t_{i}^{f},v_{i}^{f})=v_{k}^{f}(t_{i}^{f}-t_{k}^{f})-\varphi
v_{i}^{f}-\delta _{0}$ and the transversality conditions are 
\begin{equation}
\lambda _{i}^{v}(t_{i}^{f})=\eta_{i}\cdot (\frac{\partial \psi _{i}}{%
\partial v_{i}})_{t=t_{i}^{f}},H_{i}(t_{i}^{f})+\eta_{i}\cdot (\frac{%
\partial \psi _{i}}{\partial t})_{t=t_{i}^{f}}=0  \label{eq:costate_trans_2}
\end{equation}
where $\eta_i$ is the associated multiplier. (iii) When $t_{i}^{f}$ is fixed, the transversality conditions reduce to $%
\lambda _{i}^{v}(t_{i}^{f})=0$.

In addition, there also exist the state boundary conditions $%
p_{i}(t_{i}^{0})=0$, $p_{i}(t_{i}^{f})=L+S$, $v_{i}(t_{i}^{0})=v_{i}^{0}$,
given the initial time and speed $t_{i}^{0}$ and $v_{i}^{0}$.

The necessary condition for optimality is 
\begin{equation}
\frac{\partial L_{i}}{\partial u_{i}}=0=\left\{ 
\begin{array}{ll}
u_{i}(t)+\lambda _{i}^{v}(t), & \mbox{$u_{i}(t) - u_{max} <0$}~\text{and} \\ 
& \mbox{$u_{min} - u_{i}(t)<0$}, \\ 
u_{i}(t)+\lambda _{i}^{v}(t)+\mu _{i}^{a}, & \mbox{$u_{i}(t) - u_{max} =0$},
\\ 
u_{i}(t)+\lambda _{i}^{v}(t)-\mu _{i}^{b}, & \mbox{$u_{min} - u_{i}(t)=0$}.%
\end{array}%
\right.   \label{eq:optimum}
\end{equation}%
A complete solution of this problem requires that constrained and
unconstrained arcs of an optimal trajectory are pieced together to satisfy
all conditions (\ref{eq:17a}) through (\ref{eq:optimum}). This includes the
five constraints (three pure-state constraints, two control constraints) in (%
\ref{eq:17a}) through (\ref{eq:17f}). While there are many different cases
that can occur, the nature of the optimal solution rules out the possibility
of several cases. In what follows, we provide a complete analysis of the
case where no constraints are active, the case where the safety constraint $%
p_{i}(t)+\delta _{0}-p_{k}(t)\leq 0$ is the only active one, and the case
where both the state constraint $v_{i}(t)-v_{max}\leq 0$ and the control
constraint $u_{i}^{t}(t)-u_{max}\leq 0$ become active. A discussion of the
remaining cases can be found in \cite{arxiv1}.


\subsection{Unconstrained Optimal Control Analysis \label{no_active}}

For problem $P_{0}$, the terminal time is free whereas for $P_{1}$ and $P_{2}
$ the terminal time is fixed. Thus, we provide the analysis for each of
these two cases.

\subsubsection{Free Terminal Time}

\label{no_active_free} When the state and control constraints are inactive,
we have $\mu_{i}^{a}=\mu_{i}^{b}=\nu_{i}^{c}=\nu_{i}^{d}=\nu_{i}^{s}=0$. The
Lagrangian \eqref{lag} becomes $L_{i}(p,v,u,\lambda,\mu,\nu,t)=H_{i}(p,v,u,%
\lambda,t)$ and \eqref{eq:optimum} reduces to $\frac{\partial L_{i}}{%
\partial u_{i}}=u_{i}(t)+\lambda_{i}^{v}=0$, which leads to 
\begin{equation}
u_{i}(t)=-\lambda_{i}^{v}(t).   \label{optimum_1}
\end{equation}
Since $\nu_{i}^{s}=0$, \eqref{eq:lambda_p} becomes $\dot{\lambda}%
_{i}^{p}(t)=-\frac{\partial L_{i}}{\partial p_{i}}=0$ which leads to 
\begin{equation}
\lambda_{i}^{p}=a_{i}   \label{lambdaip}
\end{equation}
where $a_{i}$ is a constant. Since $\nu_{i}^{c}=\nu_{i}^{d}=0$, %
\eqref{eq:lambda_v} becomes $\dot{\lambda}_{i}^{v}(t)=-\frac{\partial L_{i}}{%
\partial v_{i}}=-\lambda_{i}^{p}$. Since $\lambda_{i}^{p}=a_{i}$, we have 
\begin{equation}
\lambda_{i}^{v}(t)=-a_{i}t-b_{i}   \label{lambdaiv}
\end{equation}
where $b_{i}$ is a constant. We can now obtain a complete analytical
solution of $P_{0}$ as follows.

The optimal trajectory for problem $P_{0}$ is given by 
\begin{align}
u_{i}^{\ast}(t) & =a_{i}t+b_{i}  \label{ui*} \\
v_{i}^{\ast}(t) & =\frac{1}{2}a_{i}t^{2}+b_{i}t+c_{i}  \label{vi*} \\
p_{i}^{\ast}(t) & =\frac{1}{6}a_{i}t^{3}+\frac{1}{2}b_{i}t^{2}+c_{i}t+d_{i} 
\label{pi*}
\end{align}
for $t\in\lbrack t_{i}^{0},t_{i}^{f^{\ast}}]$ where $a_{i}$, $b_{i}$, $c_{i}$
and $d_{i}$ are constants determined along with $t_{i}^{m^{\ast}}$ through 
\begin{subequations}
\begin{align}
\frac{1}{6}a_{i} \cdot (t_{i}^{0})^{3} + \frac{1}{2}b_{i} \cdot
(t_{i}^{0})^{2} + c_{i} t_{i}^{0} + d_{i} & = 0  \label{P0_solution:p0} \\
\frac{1}{2}a_{i} \cdot (t_{i}^{0})^{2} + b_{i} t_{i}^{0} + c_{i} & =
v_{i}^{0}  \label{P0_solution:v0} \\
\frac{1}{6}a_{i} \cdot (t_{i}^{f})^{3} + \frac{1}{2}b_{i} \cdot
(t_{i}^{f})^{2} + c_{i} t_{i}^{f} + d_{i} & = L  \label{P0_solution:pm} \\
a_{i} t_{i}^{f} + b_{i} = 0  \label{P0_solution:vm_free} \\
\gamma- \frac{1}{2}b_{i}^{2} + a_{i} c_{i} & = 0   \label{P0_solution:tran}
\end{align}
\label{P0_solution}
\end{subequations}
The optimal control in (\ref{ui*}) follows from \eqref{optimum_1} and (\ref%
{lambdaiv}). Using \eqref{ui*} in the system dynamics \eqref{eq:model2}, we
then derive (\ref{vi*}) and (\ref{pi*}). Next, (\ref{P0_solution:p0})
through (\ref{P0_solution:pm}) follow from the boundary conditions $%
p_{i}(t_{i}^{0})=0$, $v_{i}(t_{i}^{0})=v_{i}^{0}$, $p_{i}(t_{i}^{f})=L+S$
and (\ref{P0_solution:vm_free}) follows from $\lambda_{i}^{v}(t_{i}^{f})=0$
in \eqref{eq:costate_trans} and from (\ref{ui*}). The last equation follows
from $H_{i}(t_{i}^{f})=0$ in \eqref{eq:costate_trans}: 
\begin{align*}
& \gamma+\frac{1}{2}(u_{i}^{\ast}(t_{i}^{f}))^{2}+a_{i}v_{i}^{%
\ast}(t_{i}^{f})-(u_{i}^{\ast}(t_{i}^{f}))^{2} \\
& =\gamma-\frac{1}{2}(a_{i}t_{i}^{f}+b_{i})^{2}+a_{i}(\frac{1}{2}%
a_{i}(t_{i}^{f})^{2}+b_{i}t_{i}^{f}+c_{i}) \\
& =\gamma-\frac{1}{2}b_{i}^{2}+a_{i}c_{i}=0
\end{align*}
using (\ref{hamil}), (\ref{lambdaip}), (\ref{lambdaiv}), (\ref{ui*}) and (%
\ref{vi*}). 

Thus, a complete solution of $P_{0}$ boils down to solving the five
equations in (\ref{P0_solution}). A typical simulation example of this case
can be found in Section \ref{analytic} (Fig. \ref{fig:no_constraint}).

The next two results establish a basic property of the optimal control,
i.e., it is non-negative and non-increasing, and the fact that two of the
constraints in (\ref{speed_accel constraints}) cannot be active.


\begin{customlem}
For the unconstrained problem with free terminal time, the optimal control
is non-negative, i.e., $u_{i}^{\ast}(t)\geq0$, and monotonically
non-increasing \label{lem:u_nonnegative}
\end{customlem}

\emph{Proof. } Refer to \cite{arxiv1}.

\begin{customlem}
For the unconstrained problem with free terminal time, it is not possible
for constraints $v_{min}- v_{i}(t) \leq0$ and/or $u_{min} - u_{i}(t) \leq0$
to become active. \label{lem:rule_out_dec}
\end{customlem}

\emph{Proof. } Refer to \cite{arxiv1}.


\subsubsection{Constrained Terminal Time}

\label{no_active_fixed_tm} If the terminal time $t_{i}^{f}$ obtained from
solving $P_{0}$ turns out to violate \eqref{tif_lw} or \eqref{tm_up_constraint}, then, as described in Algorithm \ref{algo:optimization}, we need to solve  $P_{1}$ by setting $t_{i}^{f}$ to either the lower
bound in \eqref{tif_lw} or the upper bound in \eqref{tm_up_constraint}.
There are three subcases to consider: (i) $t_{i}^{f}$ is set to either $%
t_{i}^{L}$ or $t_{i}^{U}$, (ii) $t_{i}^{f}$ is set to a fixed value other
than $t_{i}^{L}$ and $t_{i}^{U}$, (iii) $t_{i}^{f}$ is constrained by CAV $k$%
, i.e., $t_{i}^{f}\geq t_{k}^{f}+\frac{\varphi v_{i}^{f}+\delta _{0}}{%
v_{k}^{f}}$. When $t_{i}^{f}=t_{i}^{L}$ or $t_{i}^{f}=t_{i}^{U}$, CAV $i$
simply accelerates at $u_{max}$ until reaching $v_{max}$ or decelerates at $%
u_{min}$ until reaching $v_{min}$. When $t_{i}^{f}$ is set to a fixed value
other than $t_{i}^{L}$ and $t_{i}^{U}$, the transversality conditions $%
H_{i}(t_{i}^{f})=0$ in \eqref{eq:costate_trans}, i.e., the fifth equation in
(\ref{P0_solution}), no longer holds and the solution of this problem
reduces to 
\begin{equation}
\left[ 
\begin{array}{cccc}
\frac{1}{6}(t_{i}^{0})^{3} & \frac{1}{2}(t_{i}^{0})^{2} & t_{i}^{0} & 1 \\ 
\frac{1}{2}(t_{i}^{0})^{2} & t_{i}^{0} & 1 & 0 \\ 
\frac{1}{6}(t_{i}^{f})^{3} & \frac{1}{2}(t_{i}^{f})^{2} & t_{i}^{f} & 1 \\ 
t_{i}^{f} & 1 & 0 & 0%
\end{array}%
\right] .\left[ 
\begin{array}{c}
a_{i} \\ 
b_{i} \\ 
c_{i} \\ 
d_{i}%
\end{array}%
\right] =\left[ 
\begin{array}{c}
0 \\ 
v_{i}^{0} \\ 
L+S \\ 
0%
\end{array}%
\right]   \label{fixed_term_time_solution}
\end{equation}%
which yields the four parameters $a_{i}$, $b_{i}$, $c_{i}$, $d_{i}$ from a
simple system of linear equations. A typical simulation example of this case
when \eqref{Theorem1} is violated can be found in Section \ref{analytic} (Fig. \ref{fig:no_constraint}).

When $t_{i}^{f}$ is constrained by CAV $k$, i.e., $t_{i}^{f}=t_{k}^{f}+\frac{%
\varphi v_{i}^{f}+\delta _{0}}{v_{k}^{f}}$, the transversality conditions in %
\eqref{eq:costate_trans_2} hold. We need to replace the last two equations
in \eqref{P0_solution} with the two transversality conditions in %
\eqref{eq:costate_trans_2}, i.e., $a_{i}t_{i}^{f}+b_{i}+\eta_{i}\varphi
=0$ and $\gamma +\frac{1}{2}(a_{i}t_{i}^{f}+b_{i})^{2}+\eta_{i}v_{k}^{f}=0$, where $\eta_{i}$ is the associated multiplier. In addition, we need to add the boundary
condition $t_{i}^{f}=t_{k}^{f}+\frac{\varphi v_{i}^{f}+\delta _{0}}{v_{k}^{f}%
}$ to \eqref{P0_solution}. By solving the six equations, we can obtain $a_{i}
$, $b_{i}$, $c_{i}$, $d_{i}$ along with the terminal time $t_{i}^{f}$.

With the terminal time fixed, Lemma \ref{lem:u_nonnegative} needs to be
modified as follows.

\begin{customlem}
For the unconstrained problem with fixed terminal time, the optimal control
must be either monotonically non-increasing and $u_{i}^{\ast}(t)\geq0$, or
monotonically non-decreasing and $u_{i}^{\ast}(t)\leq0.$ \label{lem:u_fixed}
\end{customlem}

\emph{Proof. } Refer to \cite{arxiv1}.

\subsection{Constrained Optimal Control Analysis \label{active}}

Checking whether the optimal solution of the unconstrained problem $P_{0}$
or $P_{1}$ violates any of the constraints (\ref{eq:17a}) through (\ref{eq:17f}) is easily accomplished
since the unconstrained optimal control \eqref{ui*} is a linear function of
time and the optimal speed is a quadratic function of time. When this
happens, we must check whether there exists a nonempty feasible control set.
One approach followed in earlier work \cite{malikopoulos2018decentralized}
is to identify the set of all initial conditions $(t_{i}^{0},v_{i}^{0})$
such that no constraint is violated over $[t_{i}^{0},t_{i}^{f}]$ or at least
some of the constraints are not violated while the rest are explicitly dealt
with through the Lagrangian in (\ref{lag}). As shown in \cite%
{Zhang2017feasibility}, it is possible to define a Feasibility Enforcement
Zone (FEZ) which precedes the CZ such that each CAV is controlled over the
FEZ so as to reach a feasible initial condition when reaching the CZ. Here,
however, we proceed differently by following a direct approach through which
we derive explicit solutions for any feasible optimal constrained
trajectory. In so doing, we can also explicitly identify when an optimal
solution is infeasible under initial conditions $(t_{i}^{0},v_{i}^{0})$.

When the optimal solution of the unconstrained problem violates a
constraint, we need to re-solve the problem by identifying an optimal
trajectory that includes unconstrained arcs pieced together with one or more
constrained arcs such that all necessary conditions for optimality are
satisfied. For a control constraint of the form $g_{i}(u_{i},t)\leq0$ as in (%
\ref{eq:17a})-(\ref{eq:17b}), the optimal control on a constrained arc can
be simply obtained by solving $g_{i}(u_{i},t)=0$. The constraints (%
\ref{eq:17c})-(\ref{eq:17e}) in our problem are pure state constraints of
the form $h_{i}(x_{i},t)\leq0$. In this case (see \cite{bryson1975optimal}),
we define the tangency constraints 
\begin{equation}
N_{i}(x_{i},t)\triangleq\lbrack h_{i}(x_{i},t)\text{ }h_{i}^{(1)}(x_{i},t)%
\text{ }\cdots\text{ }h_{i}^{(q-1)}(x_{i},t)]^{T}=0,   \label{tangent_vp}
\end{equation}
where $h_{i}^{(k)}(x_{i},t)$ is the $k$th time derivative and $q$
derivatives are taken until we obtain an expression that explicitly depends
on the control $u_{i}$ so that 
\begin{equation}
h_{i}^{(q)}(x_{i},t)=0.   \label{tangent_u}
\end{equation}
At the junction points of constrained and unconstrained arcs, the costate
and Hamiltonian trajectories may have discontinuities. This can be
determined using the following jump conditions \cite{bryson1975optimal},
where $\tau$ denotes a junction point and $\tau^{-}$,$\tau^{+}$ denote the
left-hand side and the right-hand side limits, respectively: 
\begin{equation}
\begin{aligned} \lambda_i(\tau^-) &= \lambda_i(\tau^+) + \pi_i^T
\frac{\partial N_i(x_i,t)}{\partial x_i},\\ H_i(\tau^-) &= H_i(\tau^+) -
\pi_i^T \frac{\partial N_i(x_i,t)}{\partial t}. \end{aligned} 
\label{ipm_disc}
\end{equation}
where $N_{i}(x_{i},t)$ is the $q$-dimensional vector in (\ref{tangent_vp})
and $\pi_{i}$ is a $q$-dimensional vector of constant Lagrange multipliers
satisfying $\pi_{i}^{T}N_{i}(x_{i},t)=0$ and $\pi_{i}\geq0$, $i=1,\ldots,q$.
Consequently, the optimal control $u_{i}^{\ast}(t)$ may or may not be
continuous at junction points.

In what follows, we concentrate on three cases: (i) the rear-end safety
constraint (\ref{rearend}) becomes active, (ii) the lateral collision constraint \eqref{lateral} becomes
active, (iii) both the speed constraint $v_i(t) - v_{max} \leq 0$ and the control constraint $u_i(t) - u_{max} \leq 0$ become active.

\subsubsection{Speed-dependent rear-end safety constraint $p_{k}(t)-p_{i}^{\ast}(t) \geq \protect\varphi v_i(t) + \protect\delta_0$ becomes active\label{safety_active}}

The safety constraint is the most challenging to deal with. In this case, we
have $\mu _{i}^{a}=\mu _{i}^{b}=\nu _{i}^{c}=\nu _{i}^{d}=0$. The remaining
constraints are discussed in \cite{arxiv1}. Thus, we set $%
h_{i}(p_{i},v_{i},t)=p_{i}+\varphi v_{i}+\delta _{0}-p_{k}^{\ast }(t)$ where
we observe that $p_{k}^{\ast }(t)$ is a known explicit function of time
given by the optimal position trajectory of CAV $k$ specified in (\ref%
{P0_solution}) or (\ref{fixed_term_time_solution}) since, upon arrival of
CAV $i$ at the CZ, the optimal solution of the problem associated with $k<i$
has already been fully determined. Moreover, $%
h_{i}^{(1)}(p_{i},v_{i},t)=v_{i}+\varphi u_{i}-\frac{\partial p_{k}^{\ast
}(t)}{\partial t}=v_{i}+\varphi u_{i}-v_{k}^{\ast }(t)$ where $v_{k}^{\ast
}(t)$ is also an explicit function of time in (\ref{P0_solution}) or (\ref%
{fixed_term_time_solution}).

The following result establishes the continuity property of the optimal
control when the trajectory enters a constrained arc where $p_{i}(t)+\varphi
v_i(t) + \delta_0-p_{k}^{\ast}(t)=0$.

\begin{customthm}
The optimal control $u_{i}^{\ast}(t)$ is continuous at the junction $\tau$
of the unconstrained and safety-constrained arcs, i.e., $u_{i}^{\ast}(%
\tau^{-})=u_{i}^{\ast}(\tau^{+})$. \label{safety:u_cont}
\end{customthm}

\emph{Proof. } By assumption, the rear-end safety constraint is not active
at $t_{i}^{0}$. Hence, when the safety constraint $p_{i}(t)+\varphi v_i(t) +
\delta_0-p_{k}^{\ast }(t)\leq0$ becomes active, $\tau$ is the entry time of
the constrained arc. Since $h_i^1$ explicitly depends on the control $u_i$,
we have $q=1$, and the jump conditions in \eqref{ipm_disc} become%
\begin{align*}
\lambda_{i}^{p}(\tau^{-}) & =\lambda_{i}^{p}(\tau^{+})+\pi_{i}\frac{\partial%
}{\partial p_{i}}[p_{i}+ \varphi v_i + \delta_0-p_{k}^{\ast}(t)] \\
\lambda_{i}^{v}(\tau^{-}) & =\lambda_{i}^{v}(\tau^{+})+\pi_{i}\frac{\partial%
}{\partial v_{i}}[p_{i}+ \varphi v_i + \delta_0-p_{k}^{\ast}(t)] \\
H_{i}(\tau^{-}) & =H_{i}(\tau^{+})-\pi_{i}\frac{\partial}{\partial t}[p_{i}+
\varphi v_i + \delta_0-p_{k}^{\ast}(t)]
\end{align*}
where $\frac{\partial p_{k}^{\ast}(t)}{\partial t}=v_{k}^{\ast}(t)$ and $%
\frac{\partial v_{k}^{\ast}(t)}{\partial t}=u_{k}^{\ast}(t)$ are explicit
functions of $t$ specified through (\ref{P0_solution}) or (\ref%
{fixed_term_time_solution}). We assume that $u_{k}^{\ast}(t)$, $k<i$, is
continuous in $t$ so that, if we can establish that $u_{k}^{\ast}(t)$ is
continuous, then a simple iterative argument completes the proof. The
equations above become%
\begin{align*}
\lambda_{i}^{p}(\tau^{-}) & =\lambda_{i}^{p}(\tau^{+})+\pi_{i},\text{ \ \ }%
\lambda_{i}^{v}(\tau^{-})=\lambda_{i}^{v}(\tau^{+})+\pi_{i} \varphi, \\
H_{i}(\tau^{-}) & = H_{i}(\tau^{+})+\pi_{i} v_{k}^{\ast}(t)
\end{align*}
For $t\geq\tau^{+}$, the tangency conditions (\ref{tangent_vp})-(\ref%
{tangent_u}) with $q=1$ hold:%
\begin{align*}
p_{i}(t)+\varphi v_i(t) + \delta_0-p_{k}^{\ast}(t) & = 0 \\
v_{i}(t)+ \varphi u_i(t) -v_{k}^{\ast}(t) & =0
\end{align*}
In addition, note that the position $p_{i}(t)$ and speed $v_{i}(t)$ are
continuous functions of $t$. Combining the equations above and recalling
from (\ref{hamil}) that $H_{i}(t)=\gamma+\frac{1}{2}u_{i}^{2}(t)+%
\lambda_{i}^{p}(t)v_{i}(t)+\lambda_{i}^{v}(t)u_{i}(t)$, we get%
\begin{align*}
& \gamma+\frac{1}{2}u_{i}^{2}(\tau^{-})+\lambda_{i}^{p}(\tau^{-})v_{i}(%
\tau)+\lambda_{i}^{v}(\tau^{-})u_{i}(\tau^{-}) \\
& =\gamma+\frac{1}{2}u_{i}^{2}(\tau^{+})+\lambda_{i}^{p}(\tau^{+})v_{i}(%
\tau)+\lambda_{i}^{v}(\tau^{+})u_{i}(\tau^{+}) \\
& +\pi_{i}v_{k}^{\ast}(\tau).
\end{align*}
Following from the tangency condition $v_{i}(\tau^{+}) + \varphi u_i(\tau^+)
- v_k^*(\tau^+)=0$ and the fact that $v_i(\tau^-) = v_i(\tau^+) = v_i(\tau)$%
, we have 
\begin{gather*}
\frac{1}{2}u_{i}^{2}(\tau^{-})-\frac{1}{2}u_{i}^{2}(\tau^{+})+
\pi_i[v_i(\tau) - v_k^*(\tau)] \\
+ \lambda^v(\tau^-)u_i(\tau^-) -\lambda^v(\tau^+)u_i(\tau^+) = 0
\end{gather*}
which reduces to%
\begin{align*}
&  \frac{1}{2}u_{i}^{2}(\tau^{-})-\frac{1}{2}u_{i}^{2}(\tau^{+})+
\lambda_i^v(\tau^-) (u_i(\tau^-) - u_i(\tau^+)) \\
& =[u_{i}(\tau^{-})-u_{i}(\tau^{+})][\frac{1}{2}(u_{i}(\tau^{-})+u_{i}(%
\tau^{+}))+\lambda_{i}^{v}(\tau^{-})]=0
\end{align*}
Therefore, either $u_{i}(\tau^{-})-u_{i}(\tau^{+})=0$, or $\frac{1}{2}%
[u_{i}(\tau^{-})+u_{i}(\tau^{+})]+\lambda_{i}^{v}(\tau^{-})=0$. Assuming
that $u_{i}(\tau^{-})-u_{i}(\tau^{+})\neq0$, recall that at $\tau^{-}$ the
trajectory arc is unconstrained so that (\ref{optimum_1}) holds: $%
u_{i}(\tau^{-})=-\lambda_{i}^{v}(\tau^{-})$ and it follows that $u_{i}(\tau
^{-})-u_{i}(\tau^{+})=0$. We conclude that $u_{i}(t)$ is continuous at $\tau$
and the proof is complete. 

\hfill$\blacksquare$ 

Once an optimal trajectory for CAV $i$ enters the constrained arc $%
p_{i}(t)+\varphi v_{i}(t)+\delta _{0}-p_{k}^{\ast }(t)=0$, it may remain on
this arc through the terminal time $t_{i}^{f}$ or exit it at some point $%
\tau ^{\prime }>\tau $ and follow an unconstrained arc over $[\tau ^{\prime
},t_{i}^{f}]$. This depends on whether such an exit point $\tau ^{\prime }$
is feasible on an optimal trajectory. More generally, it is possible that an
optimal trajectory consists of a sequence of alternating unconstrained and
constrained arcs whose feasibility needs to be checked. Thus, once we
establish that an optimal trajectory contains a constrained arc, there are
two cases to consider. For simplicity, let us assume that CAV $k$ is driving
within an unconstrained arc given the optimal control $u_{k}^{\ast
}(t)=a_{k}t+b_{k}$ for $t\in \lbrack t_{k}^{0},t_{k}^{f}]$ and $u_{k}^{\ast
}(t)=0$ for $t\in (t_{k}^{f},t_{i}^{f}]$, and the corresponding optimal
speed and position trajectories are $v_{k}^{\ast }(t)=\frac{1}{2}%
a_{k}t^{2}+b_{k}t+c_{k}$, $p_{k}(t)=\frac{1}{6}a_{k}t^{3}+\frac{1}{2}%
b_{k}t^{2}+c_{k}t+d_{k}$ for $t\in \lbrack t_{k}^{0},t_{k}^{f}]$, and $%
v_{k}^{\ast }(t)=v_{k}^{f}$, $p_{k}(t)=L+S+v_{k}^{f}(t-t_{k}^{f})$ for $t\in
(t_{k}^{f},t_{i}^{f}]$.

\textbf{Case 1}: No exit point from the constrained arc. In this case, CAV $i
$ remains on the constrained arc until it reaches the MZ and we have%
\begin{equation}
u_{i}^{\ast }(t)=\left\{ 
\begin{array}{cc}
a_{i}t+b_{i} & t\in \lbrack t_{i}^{0},\tau ] \\ 
a_{i}^{k}t+b_{i}^{k}+c_{e1}e^{\frac{-t}{\varphi }} & t\in (\tau ,t_{k}^{f}]
\\ 
c_{e2}e^{\frac{-t}{\varphi }} & t\in (t_{k}^{f},t_{i}^{f}]%
\end{array}%
\right.   \label{ui*_safety}
\end{equation}%
where $a_{i}^{k}=a_{k}$, $b_{i}^{k}=b_{k}-\varphi a_{i}^{k}$. CAV $i$ enters
the safety-constrained arc at $\tau $ and stays constrained until reaching
the MZ. The optimal control $u_{i}^{\ast }(t)=a_{i}^{k}t+b_{i}^{k}+c_{e1}e^{%
\frac{-t}{\varphi }}$ is derived by solving the ODE $u_{i}(t)+\varphi \dot{u}%
_{i}(t)-u_{k}^{\ast }(t)=0$ which follows from $v_{i}(t)+\varphi
u_{i}(t)-v_{k}^{\ast }(t)=0$.  Note that
for $t\in \lbrack t_{k}^{f},t_{i}^{f}]$, CAV $i$ still travels within a
safety constrained arc. Since CAV $k$ starts to cruise with $v_{k}^{f}$ at $%
t_{k}^{f}$ (Assumption \ref{ass:constant}), $a_{k}=b_{k}=0$. Hence, the
optimal form $u_{i}^{\ast }(t)=a_{i}^{k}t+b_{i}^{k}+c_{e1}e^{\frac{-t}{%
\varphi }}$ for $t\in \lbrack \tau ,t_{k}^{f}]$ reduces to $u_{i}^{\ast
}(t)=c_{e2}e^{\frac{-t}{\varphi }}$ for $t\in \lbrack t_{k}^{f},t_{i}^{f}]$.
Note that the optimal expression of CAV $i$ may vary as $u_{k}^{\ast }(t)$, $%
v_{k}^{\ast }(t)$ and $p_{k}^{\ast }(t)$ vary, which are made known to $i$
by the coordinator.

According to \eqref{eq:model2}, $v_{i}^{\ast}(t)$ is given by (\ref{vi*})
for $t\in\lbrack t_{i}^{0},\tau]$, $v_i^*(t) = \frac{1}{2}a_i^{k} t^2 +
b_i^{k}t + c_i^k - c_{e1} \varphi e^{\frac{-t}{\varphi}}$ for $%
t\in(\tau,t_{k}^{f}]$ and $v_i^*(t) = v_k^f - c_{e2} \varphi e^{\frac{-t}{%
\varphi}}$ for $t\in (t_k^f,t_{i}^{f}]$; $p_{i}^{\ast}(t)$ is given by (\ref%
{pi*}) for $t\in\lbrack t_{i}^{0},\tau]$, $p_{i}^{\ast}(t) = \frac{1}{6}%
a_i^{k} t^3 + \frac{1}{2}b_i^{k}t^2 + c_i^k t + d_i^k + c_{e1} \varphi^2 e^{%
\frac{-t}{\varphi}}$ for $t\in(\tau,t_{k}^{f}]$, and $p_{i}^{\ast}(t) = L +
S - v_k^f t_k^f - \varphi v_k^f - \delta_0 + c_{e2} \varphi^2 e^{\frac{-t}{%
\varphi}}$ for $t\in(t_k^f,t_{i}^{f}]$, where $c_i^k = c_k - \varphi b_ik$, $%
d_i^k = d_k - \varphi c_i^k - \delta_0$. The constants $a_i$, $b_i$, $c_{i}$%
, $d_{i}$ , $c_{e1}$, $c_{e2}$ along with $\tau$ and $t_i^f$ are determined
through the initial conditions, the continuity of position, speed, and
control at $\tau$ and $t_k^f$, and the terminal conditions. Simulation
examples are given in Section \ref{analytic} (Fig. \ref{fig:safety_no_exit} - \ref{fig:safety_no_exit_packed}).


\textbf{Case 2}: There exists an exit point from the constrained arc. In
this case, letting $\tau _{1}$ denote the entry point to the constrained arc
and $\tau _{2}$ the exit point, and the optimal control is given by 
\begin{equation}
u_{i}^{\ast }(t)=\left\{ 
\begin{array}{cc}
a_{i}t+b_{i} & t\in \lbrack t_{i}^{0},\tau_1 ] \\ 
a_{i}^{k}t+b_{i}^{k}+c_{e1}e^{\frac{-t}{\varphi }} & t\in (\tau_{1},\tau_{2}] \\ 
e_{i}t+r_{i} & t\in (\tau_{2},t_{i}^{f}]%
\end{array}%
\right.   \label{ui*_safety_w_exit}
\end{equation}%
For $t\in (\tau _{2},t_{i}^{f}]$, the corresponding speed and position are
given by $v_{i}^{\ast }(t)=\frac{1}{2}e_{i}t^{2}+r_{i}t+q_{i}$ and $%
p_{i}^{\ast }(t)=\frac{1}{6}e_{i}t^{3}+\frac{1}{2}r_{i}t^{2}+q_{i}t+m_{i}$.
The constants $a_{i}$, $b_{i}$, $c_{i}$, $d_{i}$ , $c_{e1}$, $e_{i}$, $r_{i}$%
, $q_{i}$, $m_{i}$, along with $\tau _{1}$, $\tau _{2}$, $t_{i}^{f}$ can be
determined through the initial conditions, the continuity of position,
speed, control at $\tau _{1}$ and $\tau _{2}$, the terminal conditions. In
terms of the terminal conditions, there are two subcases to consider: $(i)$
when the terminal time $t_{i}^{f}$ is free, and $(ii)$ when the terminal
time is fixed. When the terminal time is free, the transversality condition %
\eqref{eq:costate_trans} holds and we have $\lambda _{i}^{v}(t_{i}^{f})=0$
and $H_{i}(t_{i}^{f})=0$. In the case where the terminal time $t_{i}^{f}$ is
fixed, we simply use $\lambda _{i}^{v}(t_{i}^{f})=0$. A simulation example
is given in Section \ref{analytic} (Fig. \ref{fig:safety_w_exit}).

\begin{customrem}
Note that $a_{i}$, $b_{i}$, $c_{i}$, $d_{i}$ and $\tau_{1}$ in %
\eqref{ui*_safety_w_exit} can be determined 
independently from $e_{i}$, $r_{i}$, $q_{i}$, $m_{i}$ and $\tau_{2}$ 
if the safety constrained arc is the only constraint that becomes active.
Thus, the construction of an optimal trajectory is obtained by solving two
sub-problems and piecing the solutions together. This is an important
property because it also allows us to easily check for the existence of a
feasible solution: if $\tau_{2}<\tau_{1}$ then no feasible optimal
trajectory exists in this case.
\end{customrem}


\subsubsection{Lateral collision constraint $p_i(t_c^f) - L \leq 0$ becomes active\label{tims_active}}

When the lateral collision constraint \eqref{lateral2} becomes active, we have $t_{i}^{m}=t_{c}^{f}$ and $p_{i}(t_{c}^{f})=L$.

\begin{customthm}
When the lateral constraint \eqref{lateral2} is active at the interior-point $%
t_i^m = t_c^f$, the optimal control is continuous, i.e., $u_i^*(t_i^{m-}) =
u_i(t_i^{m+})$.\label{ucont_ipt}
\end{customthm}

\emph{Proof. } Since $N_{i}(x_{i},t_{c}^{f})=p_{i}(t_{c}^{f})-L=0$, we can
derive $\lambda _{i}^{v}(t_{i}^{m-})=\lambda _{i}^{v}(t_{i}^{m+})$ from the
jump conditions \eqref{ipm_disc}, hence, $\lambda _{i}^{v}(t)$ is continuous
at $t_{i}^{m}$. From \eqref{eq:optimum}, we know that $u_{i}^{\ast
}(t)+\lambda _{i}^{v}=0$. Since we have $u_{i}(t_{i}^{m-})+\lambda
_{i}^{v}(t_{i}^{m-})=0$ and $u_{i}(t_{i}^{m+})+\lambda _{i}^{v}(t_{i}^{m+})=0
$, we can reach the conclusion that $u_{i}(t)$ is also continuous at $%
t_{i}^{m}$.

\hfill$\blacksquare$

The optimal control for $t\in[t_i^m, t_i^f]$ can be derived in a similar way
as \eqref{ui*}, i.e., $u_i^{mz*}(t) = e_i t + r_i$, and the corresponding
speed and position are $v_i^{mz*} =\frac{1}{2}e_i t^2 + r_i t + q_i$, $%
p_i^{mz*} = \frac{1}{6}e_{i} t^3 + \frac{1}{2}r_{i} t^2 + q_{i}t + m_{i}$.
Therefore, we need four more equations for \eqref{P0_solution} to solve for $%
a_i$ through $m_i$ along with $t_i^f$ for free terminal time case (different
transversality conditions apply when $t_i^f$ is constrained), i.e., $%
p_i^*(t_i^m) = p_i^{mz*}(t_i^m)$, $v_i^*(t_i^m) = v_i^{mz*}(t_i^m)$, $%
u_i^*(t_i^m) = u_i^{mz*}(t_i^m)$ and $p_i^*(t_i^m) = L$. A simulation
example is given in Section \ref{analytic} (Fig. \ref{fig:interior_point}).


\subsubsection{Both the speed constraint $v_i(t) - v_{max} \leq 0$ and the control constraint $u_i(t) - u_{max} \leq 0$ become active\label{uv_active}}

For this case, let's consider a particular scenario where CAV $i$ will enter the arc of $u_{i}(t)-u_{max}=0$ first, and then the arc of $v_{i}(t)-v_{max}=0$. 

\begin{customthm}
The optimal trajectory cannot enter the constrained arc $v_{i}(t) - v_{max} = 0$
directly from the constrained arc $u_{i}(t) - u_{max} = 0$ if $t_i^L < t_i^f
< t_i^U$.
\end{customthm}

\emph{Proof. } First, assume that the trajectory enters the constrained arc $%
v_{i}(t)-v_{max}=0$ directly from the constrained arc $u_{i}(t)-u_{max}=0$
at $\tau $. At $\tau $, the jump conditions \eqref{ipm_disc} become 
\begin{equation}
\begin{aligned} &\lambda_i^p(\tau^-) = \lambda_i^p(\tau^+), \\\
&\lambda_i^v(\tau^-) = \lambda_i^v(\tau^+) + \pi_i(\tau), \\ &H_i(\tau^-) = H_i(\tau^+) , \\ &\pi_i(\tau)\geq 0,\pi_i(\tau)(v_i(\tau) - v_{max}) = 0. \nonumber \end{aligned}
\end{equation}%
Hence, $\lambda _{i}^{p}(t)$ and $H_{i}(t)$ are continuous at $\tau $. Since 
$v_{i}(t)$ cannot be discontinuous, and we know $u_{i}(\tau ^{-})=u_{max}$
and $u_{i}(\tau ^{+})=0$, from $H_{i}(\tau ^{-})=H_{i}(\tau ^{+})$, we have 
\begin{equation}
\begin{aligned} &\frac{1}{2}u_i(\tau^-)^2 + \lambda_i^p(\tau^-) v_i(\tau^-)
+ \lambda_i^v(\tau^-)u(\tau^-) \\ & - \frac{1}{2}u(\tau^+)^2 +
\lambda_i^p(\tau^+) v_i(\tau^+) + \lambda^v(\tau^+)u(\tau^+) \\ = &
\frac{1}{2}u_i(\tau^-)^2 + \lambda_i^v(\tau^-)u(\tau^-) - (
\frac{1}{2}u(\tau^+)^2 + \lambda^v(\tau^+)u(\tau^+)) \nonumber \end{aligned}
\end{equation}%
which reduces to 
\begin{equation}
\begin{aligned} & \frac{1}{2}u^2_{max} + \lambda^v(\tau^-)u_{max} =
u_{max}(\frac{1}{2}u_{max} + \lambda^v(\tau^-)) = 0. \nonumber \end{aligned}
\end{equation}%
Hence, we have either $u_{max}=0$ or $\frac{1}{2}u_{max}+\lambda ^{v}(\tau
^{-})=0$. If $u_{max}=0$, then CAV is not allowed to accelerate, and it is
not possible to reach $v_{max}$. If $\frac{1}{2}u_{max}+\lambda ^{v}(\tau
^{-})=0$, from \eqref{eq:optimum}, we have $\mu _{i}^{a}(\tau ^{-})=-\frac{1%
}{2}u_{max}<0$, which contradicts to $\mu _{i}^{a}(t)\geq 0$. Therefore, we
can prove that the CAV cannot enter the constrained arc $v(t)=v_{max}$
directly from the constrained arc $u(t)=u_{max}$. There exists an
unconstrained arc between the two constrained arcs. 

\hfill $\blacksquare $

\begin{customrem}
Note that if $t_i^f = t_i^L$, then CAV $i$ will simply accelerate at $u_{max}
$ until it reaches $v_{max}$. Similarly for the case when $t_i^f = t_i^U$,
CAV $i$ will decelerate at $u_{min}$ until it reaches $v_{min}$. These two
cases can be viewed as the extreme cases when the interval of the
unconstrained arc in-between reduces to zero.
\end{customrem}

Similarly to Theorem \ref{safety:u_cont}, we can also prove that $u_{i}(t)$
is continuous at both $\tau_{1}$ and $\tau_{2}$. Hence, the optimal control
is given by 
\begin{equation}
u_{i}^{\ast}(t)=\left\{ 
\begin{array}{cc}
u_{max} & t\in\lbrack t_{i}^{0},\tau_1] \\ 
a_{i}t+b_{i} & t\in(\tau_1,\tau_2] \\ 
0 & t\in(\tau_2,t_{i}^{f}]%
\end{array}
\right.   \label{ui*_uv}
\end{equation}
The coefficients can be determined through the boundary conditions (i.e.,
initial condition at $t_i^0$, terminal and transversality conditions at $%
t_i^f$) and the continuity (i.e., the continuity of position, speed and
control at $\tau_1$ and $\tau_2$). Similarly, there are also two subcases to
consider: when the terminal time is free, the transversality condition%
\eqref{eq:costate_trans} holds and we have $\lambda_i^v(t_i^f) = 0$ and $%
H_i(t_i^f) = 0$; in the case where the terminal time $t_{i}^{f}$ is fixed,
we simply use $\lambda_i^v(t_i^f) = 0$. A simulation example is given in
Section \ref{analytic} (Fig. \ref{fig:uvmax_free}).


\section{Simulation Examples\label{analytic}}

We provide several numerical examples illustrating the different cases
discussed in Section \ref{CZ}. Since the optimal solution can be obtained in decentralized fashion, with each CAV only requiring information from a subset of other CAVs, the computational time is less than 1 sec. In terms of computational complexity, we
should point out that except for the case where the complete solution is
given by the simple system of linear equations (\ref%
{fixed_term_time_solution}), solving a system of nonlinear equations
involved as in (\ref{ui*_safety}) is certainly nontrivial. A good initial
`guess' of the parameter values is extremely useful in the convergence of
the root-finding algorithm for numerical solvers. To do so, our approach is
using $u_{i}(t)$ obtained from solving problem $P_{r}$ in Algorithm \ref%
{algo:optimization}, as the initial estimation for problem $P_{r+1}$.

\subsection{Unconstrained optimal control with free terminal time.}

\begin{figure}[ptb]
\centering
\includegraphics[width= 1\columnwidth]{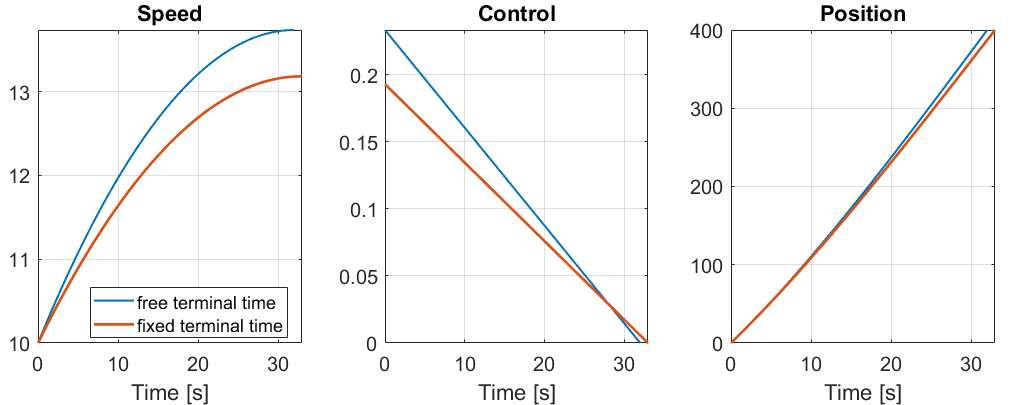} 
\caption{Unconstrained optimal trajectories with free and fixed terminal
times.}
\label{fig:no_constraint}
\end{figure}

The parameters used are: $L=370$m, $S=30$m, $\gamma=0.1$, $v_{i}^{0}=10$m/s, 
$t_{i}^{0}=0$s. The optimal terminal time is obtained as $t_{i}^{f}=32.03$s
as shown by the blue curves in Fig. \ref{fig:no_constraint}.

\textit{Unconstrained optimal control with fixed terminal time.} Assuming $%
t_{i}^{f}=32.03$s violates \eqref{tif_lw}, and we need to formulate problem $%
P_{1}$ by adding $t_{i}^{f}=33$s to $P_{0}$. The resulting optimal control,
speed, and position trajectories are shown by the red curves in Fig. \ref{fig:no_constraint}.

\subsection{Rear-end safety-constrained optimal control without exit.}

\begin{figure}[ptb]
\centering
\includegraphics[width= 1\columnwidth]{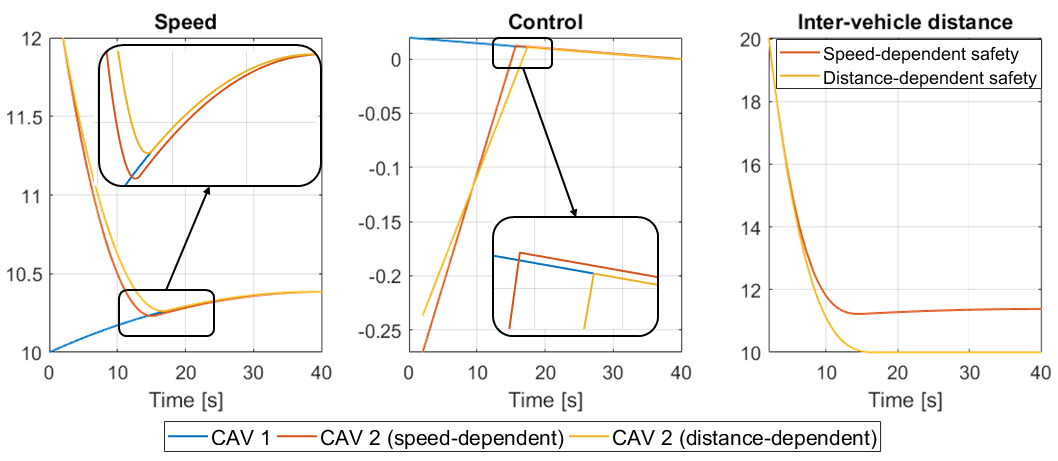}
\caption{The speed-dependent rear-end constraint $p_{i}(t)+\protect\varphi %
v_i(t) + \protect\delta_0-p_{k}(t)\leq0$ becomes active (no exit): example
\#1.}
\label{fig:safety_no_exit}
\end{figure}

\begin{figure}[ptb]
\centering
\includegraphics[width= 1\columnwidth]{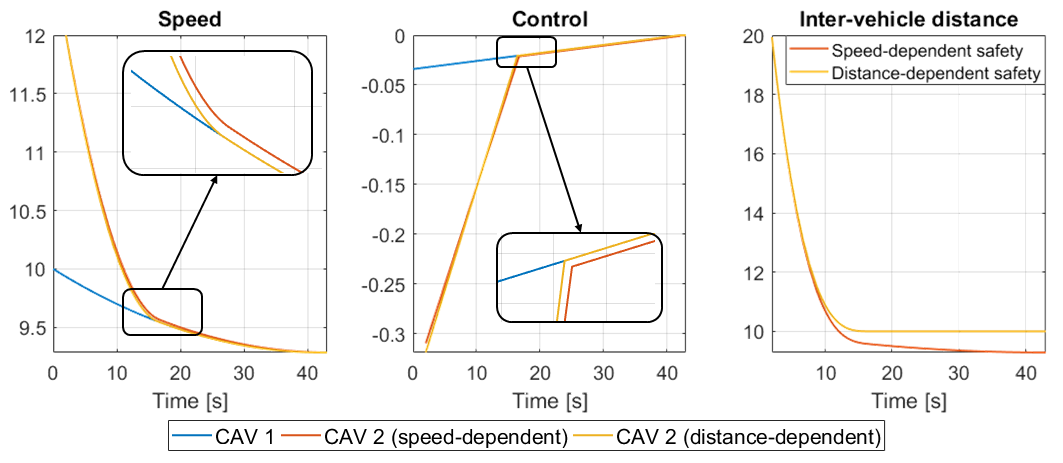}
\caption{The speed-dependent rear-end constraint $p_{i}(t)+\protect\varphi v_i(t) + \protect\delta_0-p_{k}(t)\leq0$ becomes active (no exit): example \#2.}
\label{fig:safety_no_exit_packed}
\end{figure}

Assuming CAV $k=1$ enters the CZ at $t_{k}^{0}=0$ with an initial speed $%
v_{k}^{0}=10$m/s and exits at $t_k^f=39$s, the optimal profiles are shown as
the blue curves in Fig. \ref{fig:safety_no_exit}. Then, we assume that CAV $%
i=2$ enters the CZ at $t_{i}^{0}=2$s with an initial speed $v_{i}^{0}=12$%
m/s. The coefficients for the safety constraint \eqref{rearend} is set to $%
\varphi = 1$s and $\delta_0 = 0$m. The optimal profiles for CAV $i$ is shown
as the red curves in Fig. \ref{fig:safety_no_exit}. In addition, a
comparison with the distance-dependent safety constraint as we addressed in
\cite{arxiv1} is also provided, shown as the yellow curves in Fig. \ref%
{fig:safety_no_exit}, where the minimal safe following distance is set to $10
$m. 

distance-dependent rear-end safety constraint, we provide another example by
increasing the terminal time of CAV $k$ to $t_{k}^{m}=42$s. The optimal
profiles for CAV $i$ under speed-dependent and distance-dependent constraints
are shown as the red and yellow curves respectively in Fig.
\ref{fig:safety_no_exit_packed}.

Observe that when the speed of CAV $i$ is higher than $10$m/s, the
inter-vehicle distance under the speed-dependent safety constraint increases
to ensure sufficient space between CAVs $k$ and $i$. When $v_{i}(t)$ is lower,
the required inter-vehicle distance decreases (Fig. \ref{fig:safety_no_exit_packed}), hence, vehicles can move in a more compact
manner, which improves the road utilization compared to the distance-dependent
safety constraint.

\subsection{Rear-end safety-constrained optimal control with exit.}

\begin{figure}[ptb]
\centering
\includegraphics[width= 1\columnwidth]{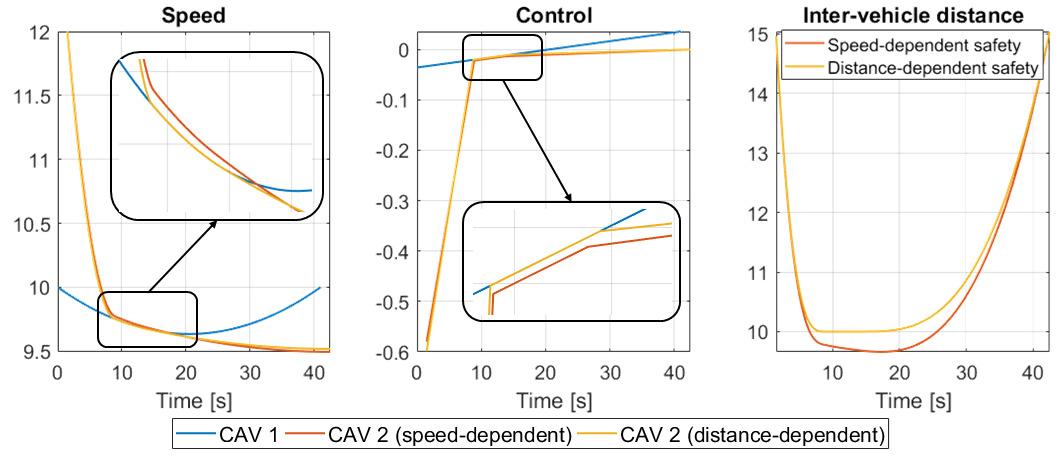}
\caption{The speed-dependent rear-end constraint $p_{i}(t) + \protect\varphi %
v_i(t) + \protect\delta_0- p_{k}(t) \leq0$ becomes active (with entry and
exit).}
\label{fig:safety_w_exit}
\end{figure}

Assuming CAV $k=1$ enters the CZ at $t_{k}^{0}=0$ with an initial speed $%
v_{k}^{0}=10$ and exits at $t_{k}^{f}=41$s with a terminal speed $%
v_{k}^{f}=10$m/s, the optimal profiles for CAV $k$ is shown as the blue
curves in Fig. \ref{fig:safety_w_exit}. Then, we assume that CAV $i=2$
enters the CZ at $t_{i}^{0}=1.5$s with an initial speed $v_{i}^{0}=12$m/s,
and the terminal time of CAV $i$ is $t_{i}^{f}=42.5$s. The optimal profiles
for CAV $i$ is shown as the red curves in Fig. \ref{fig:safety_w_exit}. 

\subsection{Lateral collision-constrained optimal control.}

\begin{figure}[ptb]
\centering
\includegraphics[width= 1\columnwidth]{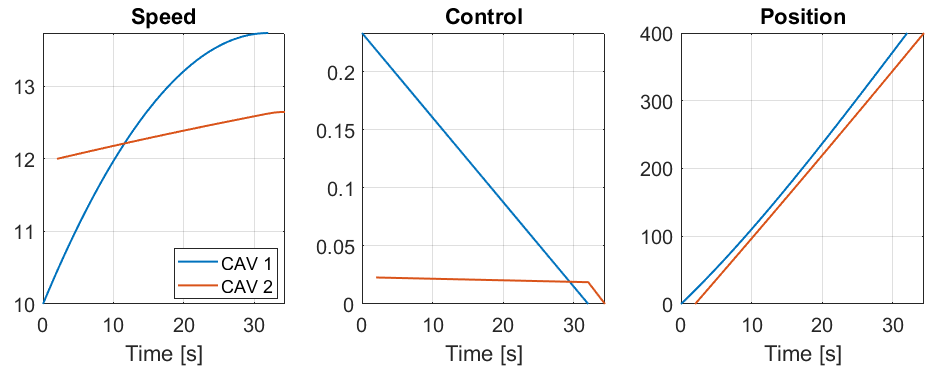}
\caption{The time-dependent lateral constraint $t_i^m \geq t_c^f$, i.e., $%
p_i(t_c^f) \leq L$ becomes active.}
\label{fig:interior_point}
\end{figure}

Assuming CAV $c=1$ enters the CZ at $t_{c}^{0}=0$ with an initial speed $%
v_{c}^{0}=10$ and exits at $t_{c}^{f}=32.027$s, the optimal profiles for CAV 
$c$ is shown as the blue curves in Fig. \ref{fig:interior_point}. Then, we
assume that CAV $i=2$ enters the CZ at $t_{i}^{0}=2$s with an initial speed $%
v_{i}^{0}=12$m/s, and the terminal time of CAV $i$ is $t_{i}^{f}=34.4$s. The
optimal profiles for CAV $i$ is shown as the red curves in Fig. \ref{fig:interior_point}. Note that CAV $c\in\mathcal{C}_{2}(t)$ and there could
be a lateral collision between them inside the MZ. Hence, CAV $i$ only enters
the MZ after CAV $c$ exits the MZ. Note that the optimal control is
continuous at $t_i^m=t_{c}^{f}=32.027$s.

\subsection{Speed-and-acceleration-constrained optimal control}

\begin{figure}[ptb]
\centering
\includegraphics[width= 1\columnwidth]{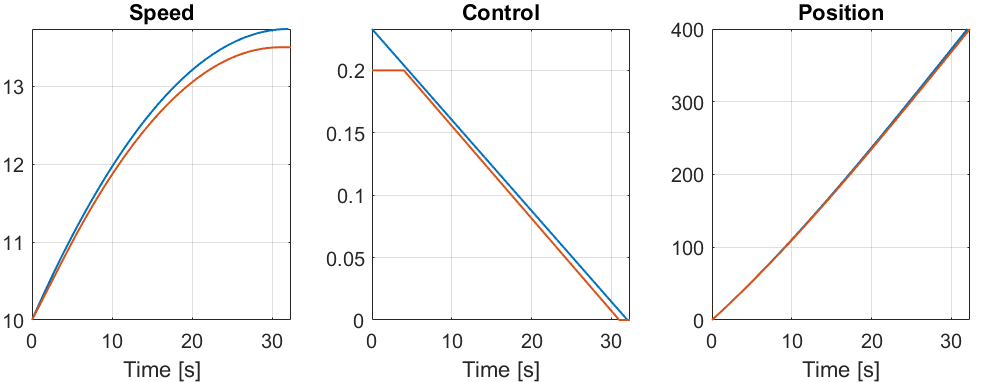}
\caption{Both the speed and control constraints $v_i(t) - v_{max} \leq 0$
and $u_i(t) - u_{max} \leq 0$ become active (free terminal time).}
\label{fig:uvmax_free}
\end{figure}

For this example, the maximum speed and acceleration are set to $v_{max} =
13.5$ and $u_{max} = 0.2$m/s$^2$, respectively. Assuming CAV $i$ enters the
CZ at $t_i^0 = 0$ with an initial speed $v_i^0 = 10$m/s. Without considering
the speed and acceleration constraints, the optimal speed exceeds the
maximum speed $v_{max}$ and the optimal control exceeds the maximum
acceleration $u_{max}$, shown as the blue curves in Fig. \ref{fig:uvmax_free}%
. Taking the constraints into consideration, the constrained optimal
trajectory consists of three arcs: one arc where CAV $i$ accelerates at $%
u_{max}$, followed by an unconstrained arc where CAV $i$ is still
accelerating but at a lower acceleration until it reaches $v_{max}$, and the
last arc where CAV $i$ cruises at $v_{max}$. The constrained optimal control
profiles are shown as the red curves in Fig. \ref{fig:uvmax_free}. Note that
the optimal control is continuous at $\tau_1 = 4.0$s, i.e., the exit point
of the control-constrained arc and $\tau_2 = 31$s, i.e., the entry point of
the speed-constrained arc. 


\section{Conclusions and Future Work}

We have extended earlier work in \cite{ZhangMalikopoulosCassandras2016} and 
\cite{malikopoulos2018decentralized} by jointly minimizing energy consumption and travel time of CAVs crossing a signal-free intersection. We include the MZ as part of the
optimal control horizon, which provides more flexibility in modeling
vehicle behavior inside the MZ. To ensure safety throughout the CZ and the
MZ, we consider a speed-dependent safety constraint, a time-dependent
lateral constraint, as well as speed and acceleration constraints, and
derive explicit solutions that possibly involve one or more of these constraints. We have also shown that the optimal solution can
still be obtained in decentralized fashion, with each CAV only requiring
information from a subset of other CAVs. This enables the on-line solution
to be obtained by on-board computation resources for each individual CAV.

Ongoing research is exploring the effect of partial CAV penetration in mixed
traffic situations where both CAVs and human-driven vehicles share the the
road \cite{Zhang2018penetration}. Future work will investigate the coupling
between multiple intersections, as well as the possibility of extending the
resequencing approach in \cite{Zhang2018sequence} to potentially improve
overall traffic throughput. 


\bibliographystyle{IEEETran}
\bibliography{CDC2019}

\end{document}